\documentclass{amsart}
\usepackage{amsmath,amscd,amsthm,amsfonts}
\usepackage[all]{xy}
\usepackage{ae}

\theoremstyle{plain}
\newtheorem{theorem}                 {Theorem}      [section]
\newtheorem{proposition}  [theorem]  {Proposition}

\newtheorem{lemma}        [theorem]  {Lemma}
\newtheorem{conjecture}        [theorem]  {Conjecture}

\theoremstyle{definition}
\newtheorem{example}      [theorem]  {Example}
\newtheorem{remark}       [theorem]  {Remark}
\newtheorem{definition}   [theorem]  {Definition}

\numberwithin{equation}{section}

\begin{document}
\baselineskip 18pt
\larger[2]

\def \theo-intro#1#2 {\vskip .25cm\noindent{\bf Theorem #1\ }{\it #2}}

\newcommand{\trace}{\operatorname{trace}}

\def \dn{\mathbb D}
\def \nn{\mathbb N}
\def \zn{\mathbb Z}
\def \qn{\mathbb Q}
\def \rn{\mathbb R}
\def \cn{\mathbb C}
\def \hn{\mathbb H}
\def \can{Ca}

\def \A{\mathcal A}
\def \B{\mathcal B}
\def \C{\mathcal C}
\def \G{\mathcal G}
\def \F{\mathcal F}
\def \H{\mathcal H}
\def \I{\mathcal I}
\def \L{\mathcal L}
\def \M{\mathcal M}
\def \N{\mathcal N}
\def \Ol{\mathcal O}
\def \P{\mathcal P}
\def \R{\mathcal R}
\def \V{\mathcal V}
\def \W{\mathcal W}

\def\Re{\mathfrak R\mathfrak e}
\def\Im{\mathfrak I\mathfrak m}
\def\Co{\mathfrak C\mathfrak o}
\def\Or{\mathfrak O\mathfrak r}

\def \ip #1#2{\langle #1,#2 \rangle}
\def \spl#1#2{( #1,#2 )}

\def \lb#1#2{[#1,#2]}

\def \pror#1{\rn P^{#1}}
\def \proc#1{\cn P^{#1}}
\def \proh#1{\hn P^{#1}}

\def \gras#1#2{G_{#1}(\cn^{#2})}

\def \g{\mathfrak{g}}
\def \k{\mathfrak{k}}
\def \m{\mathfrak{m}}
\def \p{\mathfrak{p}}
\def \q{\mathfrak{q}}
\def \r{\mathfrak{r}}
\def \un{\mathfrak{u}}

\def \GLR#1{\text{\bf GL}_{#1}(\rn)}
\def \glr#1{\mathfrak{gl}_{#1}(\rn)}
\def \GLC#1{\text{\bf GL}_{#1}(\cn)}
\def \glc#1{\mathfrak{gl}_{#1}(\cn)}
\def \GLH#1{\text{\bf GL}_{#1}(\hn)}
\def \glh#1{\mathfrak{gl}_{#1}(\hn)}
\def \GLD#1{\text{\bf GL}_{#1}(\dn)}
\def \gld#1{\mathfrak{gl}_{#1}(\dn)}

\def \SLR#1{\text{\bf SL}_{#1}(\rn)}
\def \slr#1{\mathfrak{sl}_{#1}(\rn)}
\def \SLC#1{\text{\bf SL}_{#1}(\cn)}
\def \slc#1{\mathfrak{sl}_{#1}(\cn)}

\def \O#1{\text{\bf O}(#1)}
\def \SO#1{\text{\bf SO}(#1)}
\def \so#1{\mathfrak{so}(#1)}
\def \SOO#1#2{\text{\bf SO}(#1,#2)}
\def \SOO0#1#2{\text{\bf SO}_0(#1,#2)}
\def \soo#1#2{\mathfrak{so}(#1,#2)}
\def \SOC#1{\text{\bf SO}(#1,\cn)}
\def \SOc#1{\text{\bf SO}(#1,\cn)}
\def \soc#1{\mathfrak{so}(#1,\cn)}

\def \U#1{\text{\bf U}(#1)}
\def \u#1{\mathfrak{u}(#1)}
\def \UU#1#2{\text{\bf U}(#1,#2)}
\def \uu#1#2{\mathfrak{u}(#1,#2)}
\def \SU#1{\text{\bf SU}(#1)}
\def \su#1{\mathfrak{su}(#1)}
\def \SUU#1#2{\text{\bf SU}(#1,#2)}
\def \suu#1#2{\mathfrak{su}(#1,#2)}

\def \Sp#1{\text{\bf Sp}(#1)}
\def \sp#1{\mathfrak{sp}(#1)}
\def \Spp#1#2{\text{\bf Sp}(#1,#2)}
\def \spp#1#2{\mathfrak{sp}(#1,#2)}
\def \SpC#1{\text{\bf Sp}(#1,\cn)}
\def \spc#1{\mathfrak{sp}(#1,\cn)}

\def \d#1{\mathfrak{d}(#1)}
\def \s#1{\mathfrak{s}(#1)}
\def \sym#1{\text{Sym}(\rn^{#1})}

\def \gradh#1{\text{grad}_{\H}(#1 )}
\def \gradv#1{\text{grad}_{\V}(#1 )}

\def \nab#1#2{\hbox{$\nabla$\kern -.3em\lower 1.0 ex
     \hbox{$#1$}\kern -.1 em {$#2$}}}

\allowdisplaybreaks

\title{Harmonic morphisms from the Grassmannians\\
 and their non-compact duals}

\author{Sigmundur Gudmundsson}
\author{Martin Svensson}

\thanks{The first author is a member of EDGE, Research Training Network
HPRN-CT- 2000-00101, supported by The European Human Potential
Programme.}



\address
{Mathematics, Faculty of Science, Lund University, Box 118, S-221
00 Lund, Sweden} \email{Sigmundur.Gudmundsson@math.lu.se}
\email{M.Svensson@leeds.ac.uk}

\begin{abstract}
In this paper we give a unified framework for the construction of
complex valued harmonic morphisms from the real, complex and
quaternionic Grassmannians and their non-compact duals.  This
gives a positive answer to the corresponding open existence
problem in the real and quaternionic cases.
\end{abstract}

\maketitle

\section{Introduction}

The notion of a minimal submanifold of a given ambient space is of
great importance in differential geometry.  Harmonic morphisms
$\phi:(M,g)\to(N,h)$ between Riemannian manifolds are useful tools
for the construction of such objects, see Theorem
\ref{theo:semi-B-E} below. Harmonic morphisms are solutions to
over-determined non-linear systems of partial differential
equations determined by the geometric data of the manifolds
involved.  For this reason they can be difficult to find and have
no general existence theory, not even locally.  On the contrary,
most metrics on a $3$-dimensional domain $M^3$ do not allow any
local harmonic morphisms with values in a surface $N^2$, see
\cite{Bai-Woo-1}. This makes it interesting to find geometric and
topological conditions on the manifolds $(M,g)$ and $(N,h)$,
ensuring the existence of such maps. For the general theory of
harmonic morphisms between Riemannian manifolds, we refer to the
excellent book \cite{Bai-Woo-book} and the regularly updated
on-line bibliography \cite{Gud-bib}.

For the existence of harmonic morphisms $\phi:(M,g)\to (N,h)$ it
is an advantage that the target manifold $N$ is a surface, i.e. of
dimension $2$. In this case the problem is invariant under
conformal changes of the metric on $N^2$. Therefore, at least for local
studies, the codomain can be taken to be the standard complex plane.

It is known that, in several cases when the domain $(M,g)$ is an
irreducible Riemannian symmetric space, there exist complex valued
solutions to the problem, see for example \cite{Gud-1,Sve}.
This has led the authors to the following conjecture.

\begin{conjecture}
Let $(M^m,g)$ be an irreducible Riemannian symmetric space of
dimension $m\ge 2$. For each point $p\in M$ there exists a complex
valued harmonic morphism $\phi:U\to\cn$ defined on an open
neighbourhood $U$ of $p$. If the space $(M,g)$ is of non-compact
type then the domain $U$ can be chosen to be the whole of $M$.
\end{conjecture}

It is well-known that any holomorphic map from a K\"ahler manifold
to a Riemann surface is a harmonic morphism.  This means that the
conjecture is true whenever the domain $(M,g)$ is a Hermitian
symmetric space; in particular a complex Grassmannian
$$\SU{p+q}/\text{\bf S}(\U p\times\U q)$$ or its non-compact dual
$$\SUU pq/\text{\bf S}(\U p\times\U q),$$ which can be realized as a
bounded symmetric domain in $\cn^{pq}$.

In this paper we construct explicit complex valued harmonic
morphisms defined \emph{globally} on the non-compact irreducible
Riemannian symmetric spaces $$\SOO0 pq/\SO p\times\SO q,$$ when
$p\notin \{q,q\pm 1\}$, and $$\Spp pq /\Sp p\times\Sp q$$ when
$p\neq q$. We prove the general duality Theorem \ref{theo:dual}
for complex valued harmonic morphisms from Riemannian symmetric
spaces.  This is then employed to yield locally defined solutions
from the compact real Grassmannians $$\SO{p+q}/\SO p\times\SO q,$$
with $p\notin \{q,q\pm 1\}$, and the quaternionic Grassmannians
$$\Sp{p+q}/\Sp p\times\Sp q$$ when $p\neq q$.

Throughout this article we assume that all our objects such as
manifolds,  maps etc. are smooth, i.e. in the
$C^{\infty}$-category. For our notation concerning Lie groups we
refer to the wonderful book \cite{Kna}.

\section{Harmonic Morphisms}

We are interested in complex valued harmonic morphisms from the
real, complex and quaternionic Grassmannians and their non-compact
duals. These are Riemannian manifolds, but our methods involve
harmonic morphisms from the more general semi-Riemannian
manifolds, see \cite{O-N}.

  Let $M$ and
$N$ be two manifolds of dimensions $m$ and $n$, respectively. Then
a semi-Riemannian metric $g$ on $M$ gives rise to the notion of a
Laplacian on $(M,g)$ and real-valued harmonic functions
$f:(M,g)\to\rn$. This can be generalized to the concept of a
harmonic map $\phi:(M,g)\to (N,h)$ between semi-Riemannian
manifolds being a solution to a semi-linear system of partial
differential equations, see \cite{Bai-Woo-book}.

\begin{definition}
A map $\phi:(M,g)\to (N,h)$ between semi-Riemannian manifolds is
called a {\it harmonic morphism} if, for any harmonic function $f:U\to\rn$
defined on an open subset $U$ of $N$ with $\phi^{-1}(U)$ non-empty, the
composition $f\circ\phi:\phi^{-1}(U)\to\rn$ is a harmonic function.
\end{definition}

The following characterization of harmonic morphisms between
semi-Riemannian manifolds is due to Fuglede, and generalizes the
corresponding well-known result of \cite{Fug-1,Ish} in the
Riemannian case.  For the definition of horizontal conformality we
refer to \cite{Bai-Woo-book}.

\begin{theorem}\cite{Fug-2}
  A map $\phi:(M,g)\to (N,h)$ between semi-Rie\-mannian manifolds is a
  harmonic morphism if and only if it is a horizontally (weakly)
  conformal harmonic map.
\end{theorem}

The next result generalizes the corresponding well-known theorem
of Baird and Eells in the Riemannian case, see \cite{Bai-Eel}. It
gives the theory of harmonic morphisms a strong geometric flavour
and shows that the case when $n=2$ is particularly interesting. In
that case the conditions characterizing harmonic morphisms are
then independent of conformal changes of the metric on the surface
$N^2$.  For the definition of horizontal homothety we refer to
\cite{Bai-Woo-book}.

\begin{theorem}\cite{Gud-1}\label{theo:semi-B-E}
Let $\phi:(M^m,g)\to (N^n,h)$ be a horizontally conformal
submersion from a semi-Riemannian manifold $(M^m,g)$ to a
Riemannian manifold $(N^n,h)$. If
\begin{enumerate}
\item[i.] $n=2$, then $\phi$ is harmonic if and only if $\phi$ has
minimal fibres,
\item[ii.] $n\ge 3$, then two of the following conditions imply the other:
\begin{enumerate}
\item $\phi$ is a harmonic map,
\item $\phi$ has minimal fibres,
\item $\phi$ is horizontally homothetic.
\end{enumerate}
\end{enumerate}
\end{theorem}

In what follows we are mainly interested in complex valued
functions $\phi,\psi:(M,g)\to\cn$ from semi-Riemannian manifolds.
In this situation the metric $g$ induces the complex-valued
Laplacian $\tau(\phi)$ and the gradient $\text{grad}(\phi)$ with
values in the complexified tangent bundle $T^{\cn}M$ of $M$.  We
extend the metric $g$ to be complex bilinear on $T^{\cn} M$ and
define the symmetric bilinear operator $\kappa$ by
$$\kappa(\phi,\psi)= g(\text{grad}(\phi),\text{grad}(\psi)).$$ Two
maps $\phi,\psi: M\to\cn$ are said to be {\it orthogonal} if
$$\kappa(\phi,\psi)=0.$$  The harmonicity and horizontal
conformality of $\phi:(M,g)\to\cn$ are given by the following
relations
$$\tau(\phi)=0\ \ \text{and}\ \ \kappa(\phi,\phi)=0.$$

\begin{definition}  Let $(M,g)$ be a semi-Riemannian
manifold.  A set $$\Omega=\{\phi_i:M\to\cn\ |\ i\in I\}$$ of
complex valued functions is said to be an {\it orthogonal harmonic
family} on $M$ if for all $\phi,\psi\in\Omega$
$$\tau(\phi)=0\ \ \text{and}\ \ \kappa(\phi,\psi)=0.$$
\end{definition}

\begin{remark} For a finite orthogonal harmonic family
$\{\phi_1,\dots,\phi_k\}$ on a Riemannian manifold $(M,g)$, the
map
\begin{equation*}
\Phi=(\phi_1,\dots,\phi_k):M\to\cn^k
\end{equation*}
is a {\it pseudo horizontally (weakly) conformal} map. See for
example Definition 8.2.3 and Example 8.2.6 of \cite{Bai-Woo-book}.
\end{remark}

The next result shows that the elements of an orthogonal harmonic
family can be used to produce a variety of harmonic morphisms. The
main aim of this paper is to construct such families on the
Riemannian symmetric spaces that we are dealing with.

\begin{theorem}\cite{Gud-1}\label{theo:local-sol}
Let $(M,g)$ be a semi-Riemannian manifold and
$$\Omega=\{\phi_k:M\to\cn\ |\ k=1,\dots ,n\}$$ be a finite orthogonal
harmonic family on $(M,g)$.  Let $\Phi:M\to\cn^n$ be the map given
by $\Phi=(\phi_1,\dots,\phi_n)$ and $U$ be an open subset of
$\cn^n$ containing the image $\Phi(M)$ of $\Phi$. If
$$\tilde\F=\{F_i:U\to\cn\ |\ i\in I\}$$ is a family of holomorphic
functions then $$\F=\{\psi:M\to\cn\ |\ \psi=F(\phi_1,\dots ,\phi_n
),\ F\in\tilde\F\}$$ is an orthogonal harmonic family on $M$.
\end{theorem}

\begin{proof}
The statement is a direct consequence of the fact that if
$$\psi_1=F_1(\phi_1,\dots ,\phi_n),\ \ \psi_2=F_2(\phi_1,\dots ,\phi_n)$$
are elements of $\F$ then
\begin{eqnarray*}
\tau(\psi_r)&=&\sum_{k=1}^n\frac{\partial F_r}{\partial\phi_k}\
\tau(\phi_k)+ \sum_{k,l=1}^n\frac{\partial^2
F_r}{\partial\phi_k\partial\phi_l}\
\kappa (\phi_k,\phi_l),\\
\kappa(\psi_r,\psi_s)&=&\sum_{k,l=1}^n\frac{\partial
F_r}{\partial\phi_k} \frac{\partial F_s}{\partial\phi_l}\ \kappa
(\phi_k,\phi_l).
\end{eqnarray*}
\end{proof}

\section{The Model Spaces}

In this section we introduce models for some Riemannian symmetric
spaces which are useful for our purposes. For more details, we
refer to the classical works \cite{Kob-Nom,Hel}.  Let $\dn$ be one
of the associative division algebras of the real numbers $\rn$,
complex numbers $\cn=\{x+yi\ |\ x,y\in\rn\}$ or the quaternions
$\hn=\{z+wj\ |\ z,w\in\cn\}$ of real dimension $d=1,2,4$,
respectively.  For the quaternions $\hn$ we frequently make use of
their standard representation in $\cn^{2\times 2}$ given by
$$z+wj\mapsto  \begin{pmatrix}z & w \\ -\bar w & \bar
z\end{pmatrix}.$$ By $I_n$ we denote the $n\times n$ identity
matrix and introduce the matrix
\begin{equation*}
I_{pq}=\begin{pmatrix}-I_p & 0 \\ 0 & I_q\end{pmatrix}.
\end{equation*}
Let $\dn^{(p+q)\times p}$ be the real vector space of $(p+q)\times
p$ matrices, with entries in $\dn$, equipped with either the
semi-Euclidean inner product $\spl{}{}$ given by $$\spl
XY=\Re\{\text{trace}( X^*I_{pq} Y)\}= \Re\{ -\sum_{k,l=1}^p \bar
x_{kl} y_{kl} +\sum_{k=p+1}^{p+q}\sum_{l=1}^p \bar
x_{kl}y_{kl}\}$$ or the Euclidean one $\ip{}{}$ satisfying
$$\ip XY=\Re\{\text{trace}( X^* Y)\}
=\Re\{\sum_{k=1}^{p+q}\sum_{l=1}^p \bar x_{kl}y_{kl}\}.$$
Furthermore let $\GLD p$ be the Lie group of the invertible
$\dn^{p\times p}$ matrices and $U_{pq}(\dn)$, $U^*_{pq}(\dn)$ be
the open subsets of $\dn^{(p+q)\times p}$ given by
$$U_{pq}(\dn)=\{\begin{pmatrix}X_0\\ X_1\end{pmatrix}\in\dn^{(p+q)\times p}\ |
- X_0^* X_0+ X_1^* X_1 <0\},$$
$$U^*_{pq}(\dn)=\{\begin{pmatrix}X_0\\ X_1\end{pmatrix}\in\dn^{(p+q)\times p}\
|\
X_0^*X_0+ X_1^* X_1\in\GLD p\}.$$ By the condition $- X_0^*X_0+
X_1^*X_1 <0$ we mean that for each non-zero $x\in\dn^p$ the real number
$$x^*(-X_0^*X_0+X_1^*X_1)x=-(X_0x)^*(X_0x)+(X_1x)^*(X_1x)$$ is
negative. The Lie group $\GLD p$ acts on $U_{pq}(\dn)$ and
$U^*_{pq}(\dn)$ by multiplication on the right and the quotient
spaces
$$U_{pq}(\dn)/\GLD p\ \ \text{and}\ \ U^*_{pq}(\dn)/\GLD p$$ are
differentiable manifolds of real dimension $dp q$.

Let $\Sigma_{pq}(\dn)$ and $\Sigma^*_{pq}(\dn)$ be the closed
subsets of $\dn^{(p+q)\times p}$ given by
$$\Sigma_{pq}(\dn)=\{\begin{pmatrix}X_0\\ X_1\end{pmatrix}\in\dn^{(p+q)\times p}\
|- X_0^*X_0
+X_1^*X_1=-I_p\},$$
$$\Sigma^*_{pq}(\dn)=\{\begin{pmatrix}X_0\\ X_1\end{pmatrix}\in\dn^{(p+q)\times p}\
|\
 X_0^*X_0+X_1^*X_1=I_p\}.$$   For $K_p(\rn )=\SO p$,
$K_p(\cn )=\SU p$ and $K_p(\hn )=\Sp p$ the Lie group $K_p(\dn )$
acts on $\Sigma_{pq}(\dn)$ and $\Sigma^*_{pq}(\dn)$ by
multiplication on the right and the quotient spaces
$$\Sigma_{pq}(\dn)/K_p(\dn )\ \ \text{and}\ \
\Sigma^*_{pq}(\dn)/K_p(\dn )$$ can be identified with
$U_{pq}(\dn)/\GLD p$ and $U^*_{pq}(\dn)/\GLD p$, respectively. The
metrics $\spl{}{}$ and $\ip{}{}$ on $\dn^{(p+q)\times p}$
restricted to the closed subsets $\Sigma_{pq}(\dn)$ and
$\Sigma^*_{pq}(\dn)$ induce uniquely determined Riemannian metrics
on the quotient spaces $\Sigma_{pq}(\dn)/K_p(\dn)$ and
$\Sigma^*_{pq}(\dn)/K_p(\dn )$, making the natural projections
$$\Sigma_{pq}(\dn)\to\Sigma_{pq}(\dn)/K_p(\dn )\ \ \text{and}\ \
\Sigma^*_{pq}(\dn)\to\Sigma^*_{pq}(\dn)/K_p(\dn )$$ into
Riemannian submersions.  The Riemannian manifolds obtained this
way are the real, complex and quaternionic Grassmannians and their
non-compact duals
$$\begin{array}{ccc}
\SO{p+q}/\SO p\times\SO q & & \SOO0 pq/\SO p\times\SO q\\
\SU{p+q}/\text{\bf S}(\U
p\times\U q) & & \SUU pq/\text{\bf S}(\U p\times\U q)\\
\Sp {p+q}/\Sp p\times \Sp q & &\Spp pq/\Sp p\times \Sp q.
\end{array}$$
These are well known irreducible Riemannian symmetric spaces and
the complex cases are distinguished by the fact that they carry a
natural K\"ahler structure.  The two spaces
$$\SOO0 pq/\SO p\times\SO q\ \ \text{and}\ \ \SOO0 qp/\SO q\times\SO
p,$$ obtained from each other by interchanging $p$ and $q$, are
isomorphic. The same applies to the other pairs listed above.

\begin{theorem}\label{theo:lift} Let the manifold $U_{pq}(\dn)$
be equipped with the semi-Euclidean metric $\spl{}{}$. Then the
natural projection $$\pi:U_{pq}(\dn)\to U_{pq}(\dn)/\GLD p$$ has
the following property: if $\hat\phi:U_{pq}(\dn)\to (N,h)$ is a
$\GLD p$-invariant harmonic morphism into a Riemannian manifold
$(N,h)$, then the induced map $\phi:U_{pq}(\dn)/\GLD p\to (N,h)$
is a harmonic morphism on the quotient space.
\end{theorem}

\begin{theorem}\label{theo:lift^*} Let the manifold $U_{pq}^*(\dn)$
be equipped with the Euclidean metric $\ip{}{}$. Then the natural
projection $$\pi^*:U_{pq}^*(\dn)\to U_{pq}^*(\dn)/\GLD p$$ has the
following property: if $\hat\phi^*:U^*_{pq}(\dn)\to (N,h)$ is a
$\GLD p$-invariant harmonic morphism into a Riemannian manifold,
then the induced map $\phi^*:U^*_{pq}(\dn)/\GLD p\to (N,h)$ is a
harmonic morphism on the quotient space.
\end{theorem}

We shall now prove Theorem \ref{theo:lift}. The result of Theorem
\ref{theo:lift^*} can be proved in a similar way.  This is left to
the reader as an exercise.

\begin{proof}  For simplicity we  introduce the notation
$$x_0=\begin{pmatrix}I_p\\
0\end{pmatrix}.$$ Let $G$ denote one of the groups $\SOO0 pq$,
$\SUU pq$ or $\Spp pq$, depending on the case at hand, and write
$$M=U_{pq}(\dn)/\GLD p.$$ The group $G$ acts isometrically on
$U_{pq}(\dn)$, and as $\pi$ is equivariant, it maps any fibre of
$\pi$ into a fibre of $\pi$. As the group is transitive on M, it
is transitive on the fibres of $\pi$. By $L_g$ we denote the
multiplication from the left by the element $g\in G$ and use the
same notation for this action on $U_{pq}(\dn)$, $\Sigma_{pq}(\dn)$
and on $M$. Let
$$\tilde{\pi}:\Sigma_{pq}(\dn)\to M$$ denote the restriction of
$\pi$ to $\Sigma_{pq}(\dn)$.

The fibre over $x_0$ is the $\GLD p$-orbit of $x_0$ in $U_{pq}$,
which clearly is totally geodesic. From the isometric action of
$G$, we conclude that all the fibres of $\pi$ are totally
geodesic.

Let $\hat{X}$ be a unit vector field around $\pi(x_0)$. As $\pi$
is submersive, we can find a basic vector field $X$ around $x_0$
such that $d\pi(X)=\hat{X}$. As
$\hat{X}=d\pi(X)(x_0)=d\tilde{\pi}(X)(x_0)$, $X_{x_0}$ is a unit
vector. We have
\begin{eqnarray*}
\nabla d\pi(X,X)(x_0)&=&\nabla_Xd\pi(X)(x_0)-d\pi(\nabla_XX)(x_0)\\
&=&(\nabla_{\hat{X}}\hat{X})(\pi(x_0))-d\pi(\nabla_XX)(x_0)\\
&=&(\nabla_{\hat{X}}\hat{X})(\tilde{\pi}(x_0))-d\tilde{\pi}(\nabla_XX)(x_0)\\
&=&0,
\end{eqnarray*}
as $\tilde{\pi}$ is a Riemannian submersion and $X$ is also a
basic vector field for $\tilde{\pi}$. If, on the other hand, $V$ is
a vertical vector field around $x_0$, $$\nabla
d\pi(V,V)=-d\pi(\nabla_VV)=0,$$ as $\pi$ has totally geodesic
fibres. Since $\tau(\pi)=\trace\nabla d\pi$, we have proved that
$\tau(\pi)(x_0)=0$.  The group $G$ is transitive on
$\Sigma_{pq}(\dn)$ and $L_g$ is an isometry, so we see that
$$\tau(\pi)(L_g(x_0))=\tau(\pi\circ L_g)(x_0)=\tau(L_g\circ\pi)(x_0)=0.$$
This means that $\tau(\pi)(x)=0$ for all $x\in\Sigma_{pq}(\dn)$.

Let us now assume that $\hat{\phi}=\phi\circ\pi$ is a $\GLD p
$-invariant harmonic morphism on some open subset of
$U_{pq}(\dn)$. Without loss of generality we can assume that $x_0$
is contained in the domain of $\hat\phi$. As $d\pi$ maps its
horizontal space isometrically at $x_0$ onto the tangent space
of $M$, it is clear that $\phi$ is horizontally conformal at
$\pi(x_0)$. By translating with elements in $G$, we conclude that
$\phi$ is horizontally conformal everywhere. To prove that $\phi$
is harmonic, note that
\begin{eqnarray*}
0&=&\tau(\hat{\phi})(x_0)\\
&=&d\phi(\tau(\pi)(x_0))+\trace\nabla d\phi_{\pi(x_0)}(d\pi_{x_0},d\pi_{x_0})\\
&=&\tau(\phi)(\pi(x_0)).
\end{eqnarray*}
Once again, by translating with elements in $G$, we see that
$\phi$ is harmonic.
\end{proof}

\section{The non-compact complex cases}

In this section we give a simple description of how to construct
orthogonal harmonic families on the non-compact irreducible
Hermitian symmetric spaces
$$\SUU pq/\text{\bf S}(\U p\times\U q)=U_{pq}(\cn)/\GLC p.$$
We employ Theorem \ref{theo:lift} and lift the problem into the
subset $U_{pq}(\cn)$ of $\cn^{(p+q)\times p}$ equipped with the
semi-Euclidean metric $\spl{}{}$.  This gives the following
expressions for the operators $\tau$ and $\kappa$
\begin{eqnarray*}
\tau(\phi)&=& -4\sum_{k,l=1}^p \frac{\partial^2\phi}{\partial
z_{kl}\partial\bar z_{kl}}+4\sum_{k=p+1}^{p+q}\sum_{l=1}^p
\frac{\partial^2\phi}{\partial z_{kl}\partial\bar z_{kl}},\\
\kappa(\phi,\psi)&=&-2\sum_{k,l=1}^p
        \bigl(\frac{\partial\phi}{\partial     z_{kl}}
        \frac{\partial\psi}{\partial\bar z_{kl}}
+\frac{\partial\phi}{\partial\bar z_{kl}}
        \frac{\partial\psi}{\partial z_{kl}}\bigr)\\
& &\quad +2\sum_{k=p+1}^{p+q}\sum_{l=1}^p
       \bigl(\frac{\partial\phi}{\partial     z_{kl}}
        \frac{\partial\psi}{\partial\bar z_{kl}}
       +\frac{\partial\phi}{\partial\bar z_{kl}}
        \frac{\partial\psi}{\partial z_{kl}}\bigr).
\end{eqnarray*}

\begin{proposition}\label{prop:orth-harm-fam-C}
Let $\Phi:U_{pq}(\cn)\to\cn^{q\times p}$ be the map given by
$$\Phi:\begin{pmatrix}Z_0\\ Z_1\end{pmatrix}\mapsto Z_1\cdot Z_0^{-1}.$$ Then
the
complex valued components of $\Phi$ form an orthogonal harmonic
family of $\GLC p$-invariant functions on $U_{pq}(\cn)$.
\end{proposition}

\begin{proof}
This is a direct consequence of the formulae for $\tau$ and
$\kappa$ given above and the fact that $\Phi$ is holomorphic.
\end{proof}

The components of the map $\Phi$ are $\GLC p$-invariant, and so
induce holomorphic functions on $U_{pq}(\cn)/\GLC p$ which
constitute an orthogonal harmonic family on that space.

\section{The non-compact real cases}

We shall now introduce a method for constructing orthogonal
harmonic families on the non-compact irreducible Riemannian
symmetric spaces
$$\SOO0 pq/\SO p\times\SO q=U_{pq}(\rn)/\GLR p,$$ when
$p\notin\{q,q\pm 1\}$.  It is easily seen that our method does not
work in the special cases of $p\in\{q,q\pm 1\}$.  As in the
complex case, we employ Theorem \ref{theo:lift} and lift the
problem into the set $U_{pq}(\rn)$.

Let $p,r$ be positive integers, $s=p+2r$ and on
$\rn^{(p+p+r+r)\times p}$ introduce the coordinates
$$\begin{pmatrix}A\\ B\\ W\\ \bar W\end{pmatrix}
=\begin{pmatrix}X_0-X_{1}\\ X_0+X_{1}\\  X_{2}+iX_{3}
  \\ X_{2}-iX_{3}\end{pmatrix}\quad\text{where}\quad
\begin{pmatrix}X_0\\ X_{1}\\ X_{2}\\ X_{3}\end{pmatrix}
\in\rn^{(p+p+r+r)\times p}.$$ The semi-Euclidean metric $\spl{}{}$
on $\rn^{(p+s)\times p}$ gives the following expressions for the
operators $\tau$ and $\kappa$
\begin{eqnarray*}
\tau(\phi)&=& -4\sum_{k,l=1}^p \frac{\partial^2\phi}{\partial
a_{kl}\partial b_{kl}} +4\sum_{k=1}^{r}\sum_{l=1}^p
\frac{\partial^2\phi}{\partial w_{kl}\partial\bar w_{kl}},\\
\kappa(\phi,\psi)&=&-2\sum_{k,l=1}^p
\bigl(\frac{\partial\phi}{\partial a_{kl}}
      \frac{\partial\psi}{\partial b_{kl}}
      +\frac{\partial\phi}{\partial b_{kl}}
      \frac{\partial\psi}{\partial a_{kl}}\bigr)\\
& &\quad +2\sum_{k=1}^{r}\sum_{l=1}^p
\bigl(\frac{\partial\phi}{\partial w_{kl}}
\frac{\partial\psi}{\partial\bar w_{kl}}
+\frac{\partial\phi}{\partial\bar w_{kl}}
\frac{\partial\psi}{\partial w_{kl}}\bigr).
\end{eqnarray*}

\begin{proposition}\label{prop:M-method-R}
Let $\hat M$ be an element of the complexification $\soo pr^{\cn}$
of the Lie algebra $\soo pr$ and  define the map $\hat
\Phi:\rn^{(p+s)\times p }\to\cn^{(p+r)\times p}$ by
$$\hat \Phi:X\mapsto\begin{pmatrix}A\\ W\end{pmatrix}
+ \hat M\cdot\begin{pmatrix}B\\ \bar W\end{pmatrix}.$$
Then the complex valued components of $\hat\Phi$ constitute an
orthogonal harmonic family on $\rn^{(p+s)\times p}$.
\end{proposition}

\begin{proof}  This is a simple calculation using the above
formulae for the operators $\tau$ and $\kappa$.
\end{proof}

The following result generalizes the construction of complex valued
harmonic morphisms from the odd-dimensional real hyperbolic spaces
$$\rn H^{2r-1}=\SOO0{1}{2r-1}/\SO{1}\times\SO{2r-1}$$ presented in \cite{Gud-1}.

\begin{proposition}\label{prop:orth-harm-fam-R}
Let $\Phi:U_{ps}(\rn)\to\cn^{r\times p}$ be the map given by
$$\Phi:X\mapsto W\cdot A^{-1}.  $$ Then the complex valued
components of the map $\Phi$ form an orthogonal harmonic family of
$\GLR p$-invariant functions on $U_{ps}(\rn)$.
\end{proposition}

\begin{proof}
This is a direct consequence of Proposition \ref{prop:M-method-R}
in the case when $\hat M=0$.
\end{proof}

The next result generalizes a special case of the construction
of complex valued harmonic morphisms from the even-dimensional
real hyperbolic spaces $$\rn H^{2r-2}=\SOO0{1}{2r-2}/
\SO{1}\times\SO{2r-2}$$ presented in \cite{Sve}.

\begin{proposition}\label{prop:S-method-R}
For an integer $r\ge 2$ and $M\in\soc r$ let the map
$\hat\Phi:U_{ps}(\rn ) \to\cn^{(r-1)\times p}$ be defined by
$$\hat \Phi:X\mapsto S\cdot (W + M\cdot\bar W)\cdot A^{-1}$$
where $S$ is the matrix given by
$$S=
\begin{pmatrix}
1 & 0 & \cdots & 0 & 0 & m_{1r}\\
0 & 1 & \cdots & 0 & 0 & m_{2r}\\
\vdots & \vdots &\ddots & \vdots &\vdots & \vdots\\
0 & 0 & \cdots & 1 & 0 & m_{r-2, r}\\
0 & 0 & \cdots & 0 & 1 & m_{r-1, r}\end{pmatrix}.$$ Then
$\hat\Phi$ is independent of the last row of the matrix $X$ thus
inducing a map $\Phi:U_{p,s-1}(\rn ) \to\cn^{(r-1)\times p}$. The
complex valued components of $\Phi$ constitute an orthogonal
harmonic family of $\GLR p$-invariant functions on
$U_{p,s-1}(\rn)$.
\end{proposition}

\begin{proof}
Let $t=2p+2r$ index the last row of X.  If $k=1,\dots , p$ then
\begin{eqnarray*}
& &2\frac{\partial\hat\Phi}{\partial x_{tk}}\cdot A\\
&=&S(iE_{rk}-M\cdot iE_{rk})\\
&=&\begin{pmatrix}
1 & 0 & \cdots & 0 & 0 & m_{1r}\\
0 & 1 & \cdots & 0 & 0 & m_{2r}\\
\vdots & \vdots &\ddots & \vdots &\vdots & \vdots\\
0 & 0 & \cdots & 1 & 0 & m_{r-2, r}\\
0 & 0 & \cdots & 0 & 1 & m_{r-1, r}\end{pmatrix}
\begin{pmatrix}
0 & \dots & 0 & -im_{1r} & 0 & \dots & 0 \\
\vdots &  &\vdots & \vdots &\vdots& &\vdots\\
0 & \dots & 0 & -im_{r-1,r} & 0 & \dots & 0 \\
0 & \dots & 0 & i & 0 & \dots & 0\\
\end{pmatrix}\\
&=&0.
\end{eqnarray*}
The rest follows by Proposition \ref{prop:M-method-R}.
\end{proof}

\section{The non-compact quaternionic cases}

In this section we construct orthogonal harmonic families on the
non-compact irreducible Riemannian symmetric spaces
$$\Spp pq/\Sp p\times\Sp q=U_{pq}(\hn)/\GLH p,$$ with $p\neq q$.
Our method does not work in the special cases of $p=q$.

Let $p,r$ be positive integers and set $q=p+r$.  For elements
$$Q=\begin{pmatrix}Q_0\\ Q_1\\ Q_2\end{pmatrix}\in U_{pq}(\hn)$$
we shall use the complex notation
$$Q_0=Z+Wj,\quad Q_1=X+Yj,\quad Q_2=U+Vj,\quad P=Q_0-Q_1$$ and the standard
representation of $\hn^{n\times n}$ in $\cn^{2n\times 2n}$:
$$A+Bj\mapsto\begin{pmatrix} A & B \\ -\bar B & \bar A\end{pmatrix}.$$

\begin{lemma}\label{lemm:orth-harm-fam-H}
If $p,r$ are positive integers and $q=p+r$ then
$$\Phi:U_{pq}(\hn)\to\hn^{r\times p},\quad
\Phi:Q\mapsto Q_2(Q_0-Q_1)^{-1}$$ is a $\GLH p$-invariant harmonic
map on $U_{pq}(\hn)$.
\end{lemma}

\begin{proof} Using the complex representation of $\Phi$,
$$\Phi:Q\mapsto\begin{pmatrix} U & V \\ -\bar V & \bar U\end{pmatrix}
\begin{pmatrix}Z-X & W-Y \\ \bar Y-\bar W & \bar Z-\bar X\end{pmatrix}^{-1},$$
we see that the tension field $\tau(\Phi)$ is given by the
following expression
\begin{eqnarray*}
\tau(\Phi)&=& -4\sum_{k,l=1}^p \bigl[
\frac{\partial^2\Phi}{\partial z_{kl}\partial\bar
z_{kl}}+\frac{\partial^2\Phi}{\partial w_{kl}\partial\bar
w_{kl}}\bigr]\\
& &\quad +4\sum_{k,l=1}^p\bigl[ \frac{\partial^2\Phi}{\partial
x_{kl}\partial\bar x_{kl}}
       +\frac{\partial^2\Phi}{\partial y_{kl}\partial\bar
       y_{kl}}\bigr]\\
& &\quad\quad +4\sum_{k=1}^{r}\sum_{l=1}^p\bigl[
\frac{\partial^2\Phi}{\partial u_{kl}\partial\bar u_{kl}}
       +\frac{\partial^2\Phi}{\partial v_{kl}\partial\bar
       v_{kl}}\bigr].
\end{eqnarray*}
It is obvious that the last sum vanishes. Using the following
properties for derivations $D_1,D_2$ of the matrix algebra
$\cn^{2p\times 2p}$ and an invertible element $z$ thereof,
\begin{eqnarray*}
D_1(z^{-1})&=&-z^{-1}D_1(z)z^{-1},\\
D_2(D_1(z^{-1}))&=&z^{-1}(D_2(z)z^{-1}D_1(z)-D_2(D_1(z))\\
& &\quad\quad\quad\quad\quad\quad\quad\quad
+D_1(z)z^{-1}D_2(z))z^{-1},
\end{eqnarray*}
we see that
\begin{eqnarray*}
\tau(\Phi)&=&-4Q_2P^{-1}\Bigl( \sum_{k,l=1}^p \bigl[
\frac{\partial P}{\partial z_{kl}} P^{-1} \frac{\partial
P}{\partial\bar z_{kl}}+\frac{\partial P}{\partial \bar z_{kl}}
P^{-1}\frac{\partial P}{\partial z_{kl}}\\
&&\quad\quad\quad\quad\quad\quad +\frac{\partial P}{\partial
w_{kl}} P^{-1}\frac{\partial P}{\partial\bar
w_{kl}}+\frac{\partial P}{\partial \bar w_{kl}}
P^{-1}\frac{\partial P}{\partial w_{kl}}\bigr]\\
&&\quad -\sum_{k,l=1}^p\bigl[ \frac{\partial P}{\partial x_{kl}}
P^{-1} \frac{\partial P}{\partial\bar x_{kl}}+\frac{\partial
P}{\partial \bar x_{kl}}
P^{-1}\frac{\partial P}{\partial x_{kl}}\\
&&\quad\quad\quad\quad\quad\quad  +\frac{\partial P}{\partial
y_{kl}} P^{-1}\frac{\partial P}{\partial\bar
y_{kl}}+\frac{\partial P}{\partial \bar y_{kl}}
P^{-1}\frac{\partial P}{\partial y_{kl}} \bigr]\Bigr)P^{-1}\\
&=&-4Q_2P^{-1}\Bigl( \sum_{k,l=1}^p \bigl[
\begin{pmatrix}E_{kl} & 0 \\ 0 & 0\end{pmatrix}
P^{-1}
\begin{pmatrix}0 & 0 \\ 0 & E_{kl}\end{pmatrix}\\
&&\quad +
\begin{pmatrix}0 & 0 \\ 0 & E_{kl}\end{pmatrix}
P^{-1}
\begin{pmatrix}E_{kl} & 0 \\ 0 & 0\end{pmatrix}
+\begin{pmatrix}0 & E_{kl} \\ 0 & 0\end{pmatrix} P^{-1}
\begin{pmatrix}0 & 0 \\ -E_{kl} & 0\end{pmatrix}\\
&&\quad\quad\quad\quad  +
\begin{pmatrix}0 & 0 \\ -E_{kl} & 0\end{pmatrix}
P^{-1}
\begin{pmatrix}0 & E_{kl} \\ 0 & 0\end{pmatrix}\bigr]\\
&& -\sum_{k,l=1}^p\bigl[
\begin{pmatrix}-E_{kl} & 0 \\ 0 & 0\end{pmatrix}
P^{-1}
\begin{pmatrix}0 & 0 \\ 0 & -E_{kl}\end{pmatrix}\\
&&\quad +\begin{pmatrix}0 & 0 \\ 0 & -E_{kl}\end{pmatrix} P^{-1}
\begin{pmatrix}-E_{kl} & 0 \\ 0 & 0\end{pmatrix}
+\begin{pmatrix}0 & -E_{kl} \\ 0 & 0\end{pmatrix} P^{-1}
\begin{pmatrix}0 & 0 \\ E_{kl} & 0\end{pmatrix}\\
&&\quad\quad\quad\quad  +\begin{pmatrix}0 & 0 \\
E_{kl} & 0\end{pmatrix} P^{-1}
\begin{pmatrix}0 & -E_{kl} \\ 0 & 0\end{pmatrix}\bigr]\Bigr)\\
&=&0.
\end{eqnarray*}
This shows that $\Phi$ is harmonic.
\end{proof}

The next result generalizes the construction of complex valued
harmonic morphisms from the quaternionic  hyperbolic spaces
$$\hn H^{q}=\Spp{1}{q}/\Sp{1}\times\SO{q}$$ presented in \cite{Gud-1}.

\begin{proposition}\label{prop:orth-harm-fam-H}
Let $p,r$ be positive integers, $q=p+r$ and
$\Phi:U_{pq}(\hn)\to\cn^{r\times 2p}$ the map given by
$$\Phi:Q\mapsto\begin{pmatrix} U & V \end{pmatrix}
\begin{pmatrix}Z-X & W-Y \\ \bar Y-\bar W & \bar Z-\bar X\end{pmatrix}^{-1}.$$
Then the complex valued components of $\Phi$ form an orthogonal
harmonic family of $\GLH p$-invariant functions on $U_{pq}(\hn)$.
\end{proposition}

\begin{proof}
It follows from Lemma \ref{lemm:orth-harm-fam-H} that the
components of $\Phi$ are harmonic.  The operator $\kappa$ is given
by the following equation:
\begin{eqnarray*}
& &\kappa(\phi,\psi)\\
&=& -2\sum_{k,l=1}^p \bigl( \frac{\partial\phi}{\partial
z_{kl}}\frac{\partial\psi}{\partial\bar z_{kl}}
+\frac{\partial\phi}{\partial
w_{kl}}\frac{\partial\psi}{\partial\bar w_{kl}}
+\frac{\partial\psi}{\partial
z_{kl}}\frac{\partial\phi}{\partial\bar z_{kl}}
+\frac{\partial\psi}{\partial
w_{kl}}\frac{\partial\phi}{\partial\bar w_{kl}}
\bigr)\\
& &\quad +2\sum_{k,l=1}^p\bigl( \frac{\partial\phi}{\partial
x_{kl}}\frac{\partial\psi}{\partial\bar x_{kl}}
+\frac{\partial\phi}{\partial
y_{kl}}\frac{\partial\psi}{\partial\bar y_{kl}}
+\frac{\partial\psi}{\partial
x_{kl}}\frac{\partial\phi}{\partial\bar x_{kl}}
+\frac{\partial\psi}{\partial
y_{kl}}\frac{\partial\phi}{\partial\bar y_{kl}}
\bigr)\\
& &\quad\quad +2\sum_{k=1}^{r}\sum_{l=1}^p\bigl(
\frac{\partial\phi}{\partial
u_{kl}}\frac{\partial\psi}{\partial\bar u_{kl}}
+\frac{\partial\phi}{\partial
v_{kl}}\frac{\partial\psi}{\partial\bar v_{kl}}
+\frac{\partial\psi}{\partial
u_{kl}}\frac{\partial\phi}{\partial\bar u_{kl}}
+\frac{\partial\psi}{\partial
v_{kl}}\frac{\partial\phi}{\partial\bar v_{kl}} \bigr).
\end{eqnarray*}
It is a direct consequence of the definition of $\Phi$ that
$$\frac{\partial\phi_{ij}}{\partial\bar{u}_{kl}}=0,\quad
\frac{\partial\phi_{ij}}{\partial\bar{v}_{kl}}=0.$$ Applying the
relation $D(P^{-1})=-P^{-1}D(P)P^{-1}$ we easily see that
$$\frac{\partial\phi_{ij}}{\partial z_{kl}}
\frac{\partial\phi_{rs}}{\partial\bar{z}_{kl}}
=\frac{\partial\phi_{ij}}{\partial x_{kl}}
\frac{\partial\phi_{rs}}{\partial\bar{x}_{kl}},\quad
\frac{\partial\phi_{ij}}{\partial w_{kl}}
\frac{\partial\phi_{rs}}{\partial\bar{w}_{kl}}
=\frac{\partial\phi_{ij}}{\partial y_{kl}}
\frac{\partial\phi_{rs}}{\partial\bar{y}_{kl}}.$$ This shows that
the $\GLH p$-invariant components of $\Phi$ form an orthogonal
harmonic family on $U_{pq}(\hn)$.
\end{proof}

\section{The Duality}

In this section we show how a locally defined complex valued
harmonic morphism from a Riemannian symmetric space $G/K$ of
non-compact type gives rise to a dual locally defined harmonic
morphism from its compact dual $U/K$, and vice versa.  Recall that
any harmonic morphism between real analytic Riemannian manifolds
is real analytic, see \cite{Bai-Woo-book}.

Let $W$ be an open subset of $G/K$ and $\phi:W\to\cn$ a real
analytic map. By composing $\phi$ with the projection $G\to G/K$
we obtain a real analytic $K$-invariant map $\hat\phi:\hat W\to\cn$
from some open subset $\hat W$ of $G$. Let $G^\cn$ denote the
complexification of the Lie group $G$. Then $\hat\phi$ extends
uniquely to a $K$-invariant holomorphic map $\phi^\cn:W^\cn\to\cn$
from some open subset $W^\cn$ of $G^\cn$. By restricting this map
to $U\cap W^\cn$ and factoring through the projection $U\to U/K$,
we obtain a real analytic map $\phi^*:W^*\to\cn$ from some open
subset $W^*$ of $U/K$.

\begin{theorem}\label{theo:dual}
Let $\F$ be a family of maps $\phi:W\to\cn$ and $\F^*$ be the dual
family consisting of the maps $\phi^*:W^*\to\cn$ constructed as
above.  Then $\F^*$ is an orthogonal harmonic family if and only
if $\F$ is an orthogonal harmonic family.
\end{theorem}

\begin{proof}  Let $\g=\k+\p$ be a Cartan decomposition of the Lie
algebra of $G$, where $\k$ is the Lie algebra of $K$. Furthermore
let the left-invariant vector fields $X_1,\dots,X_n\in\p$ form a
global orthonormal frame for the distribution generated by $\p$.

Let $\hat\phi$ be the lift of $\phi:W\to\cn$, via the natural
projection $\pi:G\to G/K$, defined on the open subset
$\pi^{-1}(W)$ of $G$.  We shall now assume that $\phi$ is a
harmonic morphism, i.e.
$$\tau (\hat\phi)=\sum_{k=1}^n X_k^{2}(\hat\phi)=0,\quad \kappa
(\hat\phi,\hat\phi)=\sum_{k=1}^n X_k(\hat\phi)^2=0.$$ By
construction and by the unique continuation property of real analytic
functions, the extension $\hat{\phi}^\cn$ of $\hat\phi$ satisfies
the same equations.

The Lie algebra of $U$ has the decomposition $\un=\k+i\p$ and the
left-invariant vector fields $i\bar X_1,\dots,i\bar X_n\in\ i\p$
form a global orthonormal frame for the distribution generated by
$i\p$. Let $\hat\phi^*$ be the lift of $\phi^*:W^*\to\cn$, via the
natural projection $\pi^*:U\to U/K$, defined on the open subset
$(\pi^*)^{-1}(W^*)$ of $U$.  Then
$$\tau(\hat\phi^*)=\sum_{k=1}^n(iX_k)^{2}(\hat{\phi}^*)
=-\sum_{k=1}^nX_k^{2}(\hat{\phi}^\cn)=0,$$
$$\kappa(\hat\phi^*,\hat\phi^*)=\sum_{k=1}^n(iX_k)(\hat{\phi}^*)^2
=-\sum_{k=1}^nX_k(\hat{\phi}^\cn)^2=0.$$  This shows that $\phi^*$
is a harmonic morphism.

Let us now assume that $\F$ is an orthogonal harmonic family and
that $\phi,\psi\in\F$.  Then according to Theorem
\ref{theo:local-sol} the sum $\phi+\psi$ is a harmonic morphism.
Hence the dual maps $\phi^*$, $\psi^*$ and $(\phi+\psi)^*$ are
harmonic morphisms and we have
$$\kappa ((\phi+\psi)^*,(\phi+\psi)^*)=\kappa
(\phi^*,\phi^*)+2\kappa (\phi^*,\psi^*)+\kappa
(\psi^*,\psi^*)=0.$$ Then the relations
$\kappa(\phi^*,\phi^*)=\kappa(\psi^*,\psi^*)=0$ imply that $\kappa
(\phi^*,\psi^*)=0$ in other words the harmonic dual maps $\phi^*$
and $\psi^*$ are orthogonal.  This shows that $\F^*$ is an
orthogonal harmonic family. The converse is similar.
\end{proof}

We shall apply Theorem \ref{theo:dual} to construct orthogonal
harmonic families on the compact real, complex and hyperbolic
Grassmannians.  For that purpose we now explicitly describe the
duality in the real and hyperbolic cases.

\begin{example}\label{example:real}
For the special orthogonal group $\SO{p+q}$ we have the following
Cartan decomposition
\begin{equation*}
\so{p+q}=(\so p \oplus\so q)\oplus\p,
\end{equation*}
where
\begin{equation*}
\p=\bigg\{\begin{pmatrix} 0 & X \\ -X^t & 0\end{pmatrix}
\ \bigg|\ X\in\rn^{p\times q}\bigg\}.
\end{equation*}
Let $G$ be the connected subgroup of $\SOC{p+q}$ with Lie algebra
\begin{equation*}
\g=(\so p\oplus\so q)\oplus i\p.
\end{equation*}
It is easy to see that $G$ is the identity component of the group
\begin{equation*}
\{g\in\SOC{p+q}\ |\ g^*I_{pq}g=I_{pq}\}.
\end{equation*}
Introduce the matrix
\begin{equation*}
\eta=\begin{pmatrix} -iI_p & 0 \\ 0 & I_q \end{pmatrix}.
\end{equation*}
The map $\rho:G\to\SOO0{p}{q}$
given by $$\rho(g)=\eta g\bar{\eta}$$ is an isomorphism, and we have the
following commutative diagram:
\begin{equation*}
\xymatrix{\SOC{p+q}\supset G \ar[r]^-{\rho}_-{\cong} &
\SOO0{p}{q}\ar[d]\ar[r]
& \SOO0{p}{q}/\SO{p}\times\SO{q}\ar@{=}[d]\\
 & \Sigma_{pq}(\rn) \ar[r]^-{\pi} & \Sigma_{pq}(\rn)/\SO{p}.}
\end{equation*}

Let $\phi:W\to\cn$ be a locally defined, real analytic and
$\SO{p}\times\SO{q}$-invariant map on $\SOO0{p}{q}$. Then the
composition $\phi\circ\rho$ is a $\SO{p}\times\SO{q}$-invariant
map defined locally on $G$.  This we can extend to a unique holomorphic map
$\phi^\cn:W^\cn\to\cn$ on some open subset $W^\cn$ of $\SOC{p+q}$.
The restriction $\phi^*$ of $\phi^\cn$ to $W^*=W^\cn\cap\SO{p+q}$
is a locally defined $\SO{p}\times\SO{q}$-invariant on $\SO{p+q}$.
The maps $\phi,\phi^*$ induce maps locally defined on
$\Sigma_{pq}(\rn)$ and $\Sigma_{pq}^*(\rn)$, respectively, which
we also denote by $\phi,\phi^*$. Untangling the definitions, we
see that
$$\phi^*(\binom{X}{Y})=\phi(\binom{X}{iY}),$$
where $X\in\rn^{p\times p}$ and $Y\in\rn^{q\times p}$. By employing
Theorem \ref{theo:dual}, we see that $\phi^*$ is a harmonic morphism
if and only if $\phi$ is.
\end{example}

\begin{example}\label{example:quaternionic}
For positive integers $p,q$ and $n=p+q$ we introduce the matrices
\begin{equation*}
J=\begin{pmatrix}0 & I_n \\ -I_n & 0\end{pmatrix}, \
I_{pq}=\begin{pmatrix}-I_p & 0\\ 0 & I_q\end{pmatrix}, \ K
=\begin{pmatrix}I_{pq} & 0 \\ 0 & I_{pq}\end{pmatrix}, \
\eta=\begin{pmatrix} I_n & 0 \\ 0 & I_{pq}\end{pmatrix}.
\end{equation*}
Then it is easily seen that $$\eta J\eta=KJ,\quad \eta K\eta=K.$$
Recall the following definitions:
\begin{eqnarray*}
\SpC{n}&=&\{g\in\SLC{2n}\ |\ g^tJg=J\}\\
\Sp{n}&=&\SpC{n}\cap\U{2n}=\{g\in\SU{2n}\ |\ g J=J\bar g\}\\
\Spp{p}{q}&=&\{g\in\SLC{2n}\ |\ g J=J\bar g,\ g^*Kg=K\}.
\end{eqnarray*}
Introduce the subgroup $$G=\{g\in\SpC{n}\ |\ g^*Kg=K\}$$ of
$\SpC{n}$. It is easy to see that $\Sp{p}\times\Sp{q}$ is contained in $G$
and that the map $$\rho:G\to\Spp{p}{q},\quad g\mapsto\eta g\eta$$
establishes an isomorphism between $G$ and $\Spp{p}{q}$
which preserves $\Sp{p}\times\Sp{q}$. As in the real case we have
a commutative diagram
\begin{equation*}
\xymatrix{\SpC{n}\supset G \ar[r]^-{\rho}_-{\cong} &
\Spp{p}{q}\ar[d]\ar[r]
& \Spp{p}{q}/\Sp{p}\times\Sp{q}\ar@{=}[d]\\
 & \Sigma_{pq}(\hn) \ar[r]^-{\pi} & \Sigma_{pq}(\hn)/\Sp{p}.}
\end{equation*}

Assume that $\phi:W\to\cn$ is a locally defined
$\Sp{p}\times\Sp{q}$-invariant harmonic morphism on $\Spp{p}{q}$
and let $\phi^*:W^*\to\cn$ be the $\Sp{p}\times\Sp{q}$-invariant
harmonic morphism locally defined on $\Sp{p+q}$, induced by
$\phi\circ\rho$.  Applying arguments similar to those we have used
in the real case yield the relation
\begin{equation*}
\begin{split}
\phi^*(\begin{pmatrix}
Z & W\\
X & Y\\
U & V\\
-\bar W & \bar Z\\
-\bar Y & \bar X\\
-\bar V& \bar U
\end{pmatrix}
)=\phi(\begin{pmatrix}
Z & -W\\
X & -Y\\
U & -V\\
\bar W & \bar Z\\
-\bar Y &-\bar X\\
-\bar V& -\bar U
\end{pmatrix}).
\end{split}
\end{equation*}

\end{example}

\section{The compact complex cases}

In this section we give a simple description of how to construct
orthogonal harmonic families on open subsets of the complex
Grassmannians
$$\SU {p+q}/\text{\bf S}(\U p\times\U q)=U^*_{pq}(\cn)/\GLC p.$$

The Euclidean metric $\ip{}{}$ on $\cn^{(p+q)\times p}$ gives the
following expressions for the operators $\tau$ and $\kappa$:
\begin{eqnarray*}
\tau(\phi)&=& 4\sum_{k=1}^{p+q}\sum_{l=1}^p
\frac{\partial^2\phi}{\partial z_{kl}\partial\bar z_{kl}},\\
\kappa(\phi,\psi)&=&2\sum_{k=1}^{p+q}\sum_{l=1}^p
       \bigl(\frac{\partial\phi}{\partial     z_{kl}}
        \frac{\partial\psi}{\partial\bar z_{kl}}
       +\frac{\partial\phi}{\partial\bar z_{kl}}
        \frac{\partial\psi}{\partial z_{kl}}\bigr).
\end{eqnarray*}
Let the set $V^*_{pq}(\cn )$ be defined by $$V^*_{pq}(\cn ) =\{\begin{pmatrix} Z_0\\
Z_1\end{pmatrix} \in U^*_{pq}(\cn)\ |\ \det Z_0\neq 0\}.$$

\begin{proposition}\label{prop:orth-harm-fam-C-compact}
Let $\Phi^*:V^*_{pq}(\cn)\to\cn^{q\times p}$ be the map given by
$$\Phi^*:\begin{pmatrix}Z_0\\ Z_1\end{pmatrix}\mapsto Z_1\cdot Z_0^{-1}.$$
Then the complex valued components of $\Phi^*$ constitute an
orthogonal harmonic family of $\GLC p$-invariant functions on
$V^*_{pq}(\cn)$.
\end{proposition}

\begin{proof}
This is a direct consequence of the fact that $\Phi^*$ is
holomorphic and the formulae for $\tau$ and $\kappa$ given above.
\end{proof}

\section{The compact real cases}

We shall now introduce a method for constructing orthogonal
harmonic families on open subsets of the real Grassmannians
$$\SO{p+q}/\SO p\times\SO q=U^*_{pq}(\rn)/\GLR p,$$ with
$p\notin\{q,q\pm 1\}$. Let $p,r$ be positive integers, $s=p+2r$
and introduce the complex coordinates
$$\begin{pmatrix}Z\\ \bar Z\\ W\\ \bar W \end{pmatrix}=
\begin{pmatrix}X_0+iX_{1}\\ X_0-iX_{1}\\ X_{2}+iX_{3}
  \\ X_{2}-iX_{3}\end{pmatrix}\ \ \text{where}\ \ \begin{pmatrix}X_0
  \\ X_1\\ X_2\\ X_3\end{pmatrix}\in\rn^{(p+p+r+r)\times p}.$$
The Euclidean metric $\ip{}{}$ on $\rn^{(p+s)\times p}$ gives the
following equalities for the operators $\tau$ and $\kappa$
\begin{eqnarray*}
\tau(\phi)&=& 4\sum_{k,l=1}^p \frac{\partial^2\phi}{\partial
z_{kl}\partial \bar z_{kl}} +4\sum_{k=1}^{s}\sum_{l=1}^p
\frac{\partial^2\phi}{\partial w_{kl}\partial\bar w_{kl}},\\
\kappa(\phi,\psi)&=&2\sum_{k,l=1}^p
\bigl(\frac{\partial\phi}{\partial z_{kl}}
      \frac{\partial\psi}{\partial \bar z_{kl}}
      +\frac{\partial\phi}{\partial \bar z_{kl}}
      \frac{\partial\psi}{\partial z_{kl}}\bigr)\\
& &\quad +2\sum_{k=1}^{s}\sum_{l=1}^p
\bigl(\frac{\partial\phi}{\partial w_{kl}}
\frac{\partial\psi}{\partial\bar w_{kl}}
+\frac{\partial\phi}{\partial\bar w_{kl}}
\frac{\partial\psi}{\partial w_{kl}}\bigr).
\end{eqnarray*}

\begin{proposition}\label{prop:M-method-R-compact}
Let the matrix $\hat M$ be an element of the complex Lie algebra
$\soc {p+r}$ and define the map $\hat \Phi^*:\rn^{(p+s)\times p
}\to\cn^{(p+r)\times p}$ by
$$\hat \Phi^*:X\mapsto\begin{pmatrix}Z\\ W\end{pmatrix}
+ \hat M\cdot\begin{pmatrix}\bar Z\\ \bar W\end{pmatrix}.$$
Then the complex valued components of $\hat\Phi^*$ constitute an
orthogonal harmonic family on $\rn^{(p+s)\times p}$.
\end{proposition}

\begin{proof}  This is a simple calculation using the above
formulae for the operators $\tau$ and $\kappa$.
\end{proof}

The following result generalizes the construction of complex valued
harmonic morphisms from the odd-dimensional real projective spaces
$$\rn P^{2r-1}=\SO{2r}/\SO{1}\times\SO{2r-1}$$ presented in \cite{Gud-1}.
For this we define the set $V^*_{pq}(\rn )$ by the
formula $$V^*_{pq}(\rn )=\{X\in U^*_{pq}(\rn)\ |\ \det Z\neq
0\}.$$

\begin{proposition}\label{prop:orth-harm-fam-R-compact}
Let $\Phi^*:V^*_{ps}(\rn)\to\cn^{r\times p}$ be the map defined by
$$\Phi^*:X\mapsto W\cdot Z^{-1}.$$  Then the complex valued
components of $\Phi^*$ form an orthogonal harmonic family of $\GLR
p$-invariant functions on $V^*_{ps}(\rn)$.
\end{proposition}

\begin{proof}
This is a direct consequence of Proposition
\ref{prop:M-method-R-compact} in the case when $\hat M=0$.
\end{proof}

\begin{proposition}\label{prop:S-method-R-compact}
For an integer $r\ge 2$ and $M\in\soc r$ let the map
$\hat\Phi^*:V^*_{ps}(\rn ) \to\cn^{(r-1)\times p}$ be defined by
$$\hat\Phi^*:X\mapsto S\cdot (W + M\cdot\bar W)\cdot Z^{-1}$$
where $S$ is the matrix given by
$$S=
\begin{pmatrix}
1 & 0 & \cdots & 0 & 0 & m_{1r}\\
0 & 1 & \cdots & 0 & 0 & m_{2r}\\
\vdots & \vdots &\ddots & \vdots &\vdots & \vdots\\
0 & 0 & \cdots & 1 & 0 & m_{r-2, r}\\
0 & 0 & \cdots & 0 & 1 & m_{r-1, r}\end{pmatrix}.$$ Then
$\hat\Phi^*$ is independent of the last row of the matrix $X$,
thus inducing a map $\Phi^*:V^*_{p,s-1}(\rn ) \to\cn^{(r-1)\times
p}$. The complex valued components of $\Phi^*$ form an orthogonal
harmonic family of $\GLR p$-invariant functions on
$V^*_{p,s-1}(\rn)$.
\end{proposition}

\begin{proof}
Let $t=2p+2r$ index the last row of X.  If $k=1,\dots , p$ then
\begin{eqnarray*}
& &2\frac{\partial\hat\Phi^*}{\partial x_{tk}}\cdot Z\\
&=&S(iE_{rk}-M\cdot iE_{rk})\\
&=&\begin{pmatrix}
1 & 0 & \cdots & 0 & 0 & m_{1r}\\
0 & 1 & \cdots & 0 & 0 & m_{2r}\\
\vdots & \vdots &\ddots & \vdots &\vdots & \vdots\\
0 & 0 & \cdots & 1 & 0 & m_{r-2, r}\\
0 & 0 & \cdots & 0 & 1 & m_{r-1, r}\end{pmatrix}
\begin{pmatrix}
0 & \dots & 0 & -im_{1r} & 0 & \dots & 0 \\
\vdots &  &\vdots & \vdots &\vdots& &\vdots\\
0 & \dots & 0 & -im_{r-1,r} & 0 & \dots & 0 \\
0 & \dots & 0 & i & 0 & \dots & 0\\
\end{pmatrix}\\
&=&0.
\end{eqnarray*}
The rest follows by Proposition \ref{prop:M-method-R-compact}.
\end{proof}

\section{The compact quaternionic cases}

In this section we construct orthogonal harmonic families on open
subsets of the quaternionic Grassmannians
$$\Sp {p+q}/\Sp p\times\Sp q=U^*_{pq}(\hn)/\GLH p,$$ with $p\neq q$.

Let $p,r$ be positive integers and set $q=p+r$.  As before, we use
the notation
$$\begin{pmatrix} Q_0\\ Q_1\\ Q_2\end{pmatrix}=
\begin{pmatrix} Z+Wj\\ X+Yj\\ U+Vj\end{pmatrix}
\ \ \text{with}\ \
\begin{pmatrix}Q_0\\ Q_1\\ Q_2\end{pmatrix}\in U^*_{pq}(\hn).$$
Define the set $$V^*_{pq}(\hn)
=\{Q\in U^*_{pq}(\hn)\ |\ \det
\begin{pmatrix}Z-X & Y-W \\ \bar Y+\bar W & \bar
Z+\bar X\end{pmatrix}\neq 0\}.$$

\begin{proposition}\label{prop:orth-harm-fam-H-compact}
Let $p,r$ be positive integers, $q=p+r$ and let
$\Phi^*:V^*_{pq}(\hn)\to\cn^{2r\times 2p}$ be the map given by
$$\Phi^*:Q\mapsto\begin{pmatrix} U & -V \end{pmatrix}
\begin{pmatrix}Z-X & Y-W \\ \bar Y+\bar W & \bar Z+\bar X\end{pmatrix}^{-1}.$$
Then the complex valued components of $\Phi^*$ constitute an orthogonal
harmonic family of $\GLH p$-invariant functions on
$V^*_{pq}(\hn)$.
\end{proposition}

\begin{proof}
The statement follows from Theorem \ref{theo:dual} combined with
Proposition \ref{prop:orth-harm-fam-H} and Example
\ref{example:quaternionic}.
\end{proof}

\end{document}